\documentclass[11pt]{amsart}

\usepackage{float}
\usepackage[margin=1in]{geometry}
\usepackage{amsmath,amsthm,amssymb}
\usepackage{bm}
\usepackage{xcolor}
\usepackage{hyperref}
\usepackage{cleveref}
\usepackage{mathrsfs}
\usepackage{graphicx}

\usepackage[
backend=biber,
style=alphabetic
]{biblatex}
\numberwithin{equation}{section}

\newcommand{\Z}{\mathbb{Z}}
\newcommand{\R}{\mathbb{R}}
\newcommand{\Q}{\mathbb{Q}}

\newcommand{\p}{\mathbb{P}}

\newtheorem{theorem}{Theorem}[section]
\newtheorem{lemma}[theorem]{Lemma}

\newtheorem{corollary}[theorem]{Corollary}

\newtheorem*{remark}{Remark}

\crefname{lemma}{lemma}{lemmas}
\crefname{proposition}{proposition}{propositions}
\crefname{theorem}{theorem}{theorems}
\crefname{corollary}{corollary}{corollaries}

\begin{document}

\title{Average density of rational points on quartic surfaces}
\author{Anna Zanoli}
\address{Institute of Science and Technology Austria (ISTA),
Am Campus 1,
3400 Klosterneuburg, Austria.}
\email{anna.zanoli@ista.ac.at}
\date{}

\begin{abstract}
We establish an average asymptotic for the number of rational points of bounded height $B$ on quartic surfaces in $\mathbb{P}^3$ of height bounded by $A$, refining an argument of Le Boudec and using techniques of Browning--Le Boudec--Sawin. We extend the range of the asymptotic from $B=o(A)$ to $B=O(A^2)$. As an application, we show that almost all quartic surfaces of height $H$ have no rational point of height $o(H^{5/2})$.
\end{abstract}
\maketitle

\tableofcontents

\section{Introduction}

Let $\mathbb{V}$ denote the family of quartic hypersurfaces in
$\p^3$ defined over $\Q$. The number of monomials of degree $4$ in $4$ variables is
$\binom{3+4}{4}=35$. Ordering the monomials lexicographically, we can therefore identify
\[
    \mathbb{V}=\p^{34}(\Q).
\]
A generic element of this family is a smooth quartic surface. By the
adjunction formula its canonical line bundle is trivial, and hence it is a
K3 surface.

For $V\in\mathbb{V}$, let
$\bm{a}_V\in\Z^{35}_{\mathrm{prim}}$ be either of the two primitive coefficient
vectors representing $V$, and let $\|\cdot\|$ denote the Euclidean norm. Define $H(V)=\|\bm{a}_V\|$. For $A\geq1$, set
\[
    \mathcal{F}_A
    :=
    \left\{V\in\mathbb{V}:H(V)\leq A\right\}.
\]
Similarly, for $x\in\p^3(\Q)$ we choose a primitive representative
$\bm{x}\in\Z^4_{\mathrm{prim}}$ and define $H(x)=\|\bm{x}\|$. Let $V_F$ denote the hypersurface $\{F(\bm{x})=0\}$ for
$F\in\Z[x_0,x_1,x_2,x_3]$ homogeneous of degree $4$. We write
\[
    N_F(B)
    :=
    \#\left\{x\in V_F(\Q):H(x)\leq B\right\}.
\]
Our object of study is the first moment
\[
    \mathbb{E}_A(N_F(B))
    :=
    \frac{1}{\#\mathcal{F}_A}
    \sum_{V_F\in\mathcal{F}_A}N_F(B).
\]

Related statistical questions for random Fano hypersurfaces have been
studied by Le Boudec \cite{leboudec2020height} and by
Browning--Le Boudec--Sawin \cite{brbosaw}. In particular, Le Boudec proves an
average version of Manin's conjecture for Fano hypersurfaces of degree $d$ in $\p^{n}$ with $n\geq d$
\cite[Theorem~3]{leboudec2020height}, and remarks that the same method extends
to hypersurfaces of intermediate type, that is, to the case $n=d-1$, upon
replacing the anticanonical height
$H(\bm{x})=\|\bm{x}\|^{n+1-d}$ by the naive height
$H(\bm{x})=\|\bm{x}\|$. The case $(d,n)=(4,3)$ considered here is precisely
the surface instance of this intermediate regime.

In this boundary case, the anticanonical line bundle becomes
trivial rather than ample and the change of geometry is reflected in the average
growth of rational points. In place of the linear growth in $B$ found on average for
generic Fano hypersurfaces, one obtains logarithmic growth
\cite[Introduction, after Theorem~3]{leboudec2020height}. More precisely,
specializing Le Boudec's calculation to $(d,n)=(4,3)$ gives, in the range
$A\geq B$,
\begin{equation}\label{eq:leboudec-intro}
    \mathbb{E}_A(N_F(B))
    =
    c_4\frac{\log B}{A}
    \left(
        1+
        O\left(
            \frac BA+\frac1{\log B}
        \right)
    \right).
\end{equation}
As Le Boudec points out, this estimate is in agreement with the unpublished
observation of van Luijk that the number of rational points of height bounded
by $B$ on a suitable K3 surface should behave asymptotically like $\log B$ to
the power of the rank of the Picard group, which for almost all elements of
$\mathbb{V}$ is equal to $1$ (this has been known since Noether and a modern proof was given by Deligne, see \cite[Theorem~XIX.1.2]{DeligneKatz1973}).

The aim of this work is to improve the error term in
\eqref{eq:leboudec-intro} and to enlarge the range of $B$ relative to $A$. Let $C_m$ be the volume of the unit ball in $m$-dimensional Euclidean space, $\nu:\R^4\to\R^{35}$ the Veronese map and $d\sigma(\bm{t})$ the standard
surface measure on $S^3$ induced by the Euclidean metric. Our main quantitative estimate is the following.

\begin{theorem}\label{thm:first-moment}
For $A,B\geq2$, we have
\[ \mathbb{E}_A(N_F(B)) = c_4\frac{\log B}{A} + O\left( \frac1A  + \frac{B^2}{A^5} + \frac{B^5}{A^{13}} \right), \]
where
\[
    c_4
    :=
    \frac12
    \frac{C_{34}\zeta(35)}
    {C_{35}\zeta(34)\zeta(4)}
    \int_{S^3}
    \frac{d\sigma(\bm{t})}{\|\nu(\bm{t})\|}.
\]
\end{theorem}

\begin{corollary}\label{cor:first-moment} Let $A,B\to\infty$. If $ B^2=o\left(A^4\log B\right), $ then \[ \mathbb{E}_A(N_F(B)) \sim c_4\frac{\log B}{A}. \] In particular, the asymptotic holds whenever $ B=O(A^2). $ \end{corollary}

Note that there is no hope to go beyond the quadratic range. To see this, fix the rational line $L:\{x_2=x_3=0\}\subset \p^3$.  The restriction of a quartic form to $L$ is a binary
quartic and has five coefficients.  Thus the forms containing $L$ form a
subfamily of $\mathcal{F}_A$ of codimension five, with density $\asymp A^{-5}$. Meanwhile
\[
\#\{x\in L(\Q):H(x)\leq B\}\asymp B^2.
\]
We deduce that this single subfamily contributes
\begin{equation*}
\gg\frac{B^2}{A^5}
\end{equation*}
to the complete family average.  This exceeds the expected average main term $\frac{\log B}{A}$ if
$B^2\gg A^4\log B$. Thus the range for $A,B$ in \Cref{cor:first-moment} is optimal.

The proof refines Le Boudec's lattice-point argument by exploiting additional information on the successive minima of a certain lattice $\Lambda_{\bm{x}}$. Revisiting the construction underlying \cite[Lemma~3.15]{brbosaw}, we obtain, in addition to the bound for the largest successive minimum, a useful estimate for smaller successive minima. We divide the vectors into a good range, where the ball of radius $A$ contains a full set of independent lattice vectors and Le Boudec's asymptotic estimate applies, and a bad range. In the latter, we combine the refined successive-minima bounds with results and techniques from \cite{brbosaw}. This gives a substantially sharper estimate for the bad contribution and extends the asymptotic to the quadratic scale $B=O(A^2)$.

\Cref{thm:first-moment} can be used to deduce information about the minimal height of rational points on random quartic surfaces. More precisely, for $V\in\mathbb{V}$, define
\[
M(V)=
\begin{cases}
\displaystyle \min_{x\in V(\mathbb Q)} H(x),
    & V(\mathbb Q)\neq\varnothing,\\[1ex]
\infty,
    & V(\mathbb Q)=\varnothing.
\end{cases}
\]
We will prove the following.
\begin{theorem}\label{thm:smallest-point}
Let $\psi:\mathbb R_{>0}\to\mathbb R_{>0}$ satisfy
$\psi(u)=o\left(u^{\frac{5}{2}}\right)$ as $u\to\infty.$ Then
\[
\lim_{A\to\infty}
\frac{
\#\{V\in\mathcal F_A:M(V)>\psi(H(V))\}
}{
\#\mathcal F_A
}
=1.
\]
\end{theorem}
The exponent $\frac{5}{2}$ should not be regarded as a prediction for
the natural scale of the smallest rational point. It arises from the
first-moment method, since the contribution $\frac{B^2}{A^5}$ of the surfaces containing a rational line becomes of order one when $B\asymp A^{\frac{5}{2}}$, but the main term itself suggests a rather different scale. Formally, if
the behaviour $\mathbb E_A(N_F(B))
    \sim
    c_4\frac{\log B}{A}$ were to persist in a much larger range, then the expected number of
rational points of height at most $B$ would become of order one when $\log B\asymp A$,
that is, for $B$ exponential in $A$. 

\medskip

In Section~\ref{sec:geometry-numbers}, we collect the geometry-of-numbers
estimates used throughout the argument, including the relevant results of
Le Boudec and Browning--Le Boudec--Sawin. In
Section~\ref{sec:first-moment}, we apply these tools to the first moment:
we first treat the good contribution, then estimate the bad
contribution, and finally combine the two estimates to prove
Theorem~\ref{thm:first-moment}. In \Cref{sec:minimalheight}, we apply \Cref{thm:first-moment} to the minimal height of a rational point on a quartic surface $V$, proving \Cref{thm:smallest-point}.

\subsection*{Acknowledgements}
The author would like to thank her advisor, Tim Browning, for suggesting
this problem, for the many helpful discussions, and for his encouragement and
support throughout the project.

\section{Geometry of numbers}\label{sec:geometry-numbers}

The aim of this section is to collect the geometry-of-numbers results
needed for the proof of our main theorem and to fix the notation that
will be used throughout.

For any $N\geq1$ and vectors $\bm{\alpha},\bm{\beta}\in\R^N$, we denote by
$\bm{\alpha}\cdot\bm{\beta}$ their scalar product. We also recall that
$\|\cdot\|$ denotes the Euclidean norm in $\R^N$. We write
\[
    \mathcal{B}_m(T)=\{\bm{y}\in\R^m:\|\bm{y}\|\leq T\},
    \qquad
    C_m=\operatorname{vol}(\mathcal{B}_m(1)).
\]
Given an integral lattice $\Lambda\subset\Z^N$, write
\[
    \mathcal{N}^*(\Lambda,\mathcal{B}_N(T))
    :=
    \#\left(
        \Lambda
        \cap
        \Z^N_{\mathrm{prim}}
        \cap
        \mathcal{B}_N(T)
    \right).
\]
The following estimate builds upon work of Schmidt
\cite[Lemma~2]{Schmidt1968}.

\begin{lemma}[{\cite[Lemma~3]{leboudec2020height}}]
\label{lemma:primitive-count}
Let $N\geq3$ and $R\in\{3,\ldots,N\}$. Let
$\Lambda\subset\R^N$ be an integral and primitive lattice of rank $R$,
and suppose that $\mathcal{B}_N(Y)$ contains $R$ linearly independent
vectors of $\Lambda$. If $T\geq Y\geq1$, then
\[
    \mathcal{N}^*(\Lambda,\mathcal{B}_N(T))
    =
    \frac{C_R}{\zeta(R)}
    \frac{T^R}{\det(\Lambda)}
    \left(
        1+
        O\left(\frac YT\right)
    \right)
    +
    O(T\log Y),
\]
where the implied constants depend at most on $R$.
\end{lemma}
Note that applying this with    $\Lambda=\Z^{35}$, $R=N=35$, $Y=1$, directly gives
\begin{equation}\label{eq:family-size}
    \#\mathcal{F}_A
    =
    \frac{C_{35}}{2\zeta(35)}
    A^{35}
    \left(
        1+
        O\left(\frac1A\right)
    \right).
\end{equation}

Now let $\nu:\R^4\longrightarrow\R^{35}$ be the degree $4$ Veronese map, obtained by listing all monomials of
degree $4$ in lexicographic order,
\[
    \nu(\bm{x})=(x_0^4,x_0^3x_1,\ldots,x_3^4).
\]
For $\bm{x}\in\Z^4_{\mathrm{prim}}$, let
\[
    \Lambda_{\bm{x}}
    =
    \left\{
        \bm{c}\in\Z^{35}:
        \bm{c}\cdot\nu(\bm{x})=0
    \right\}.
\]
Since $\bm{x}$ is primitive, $\nu(\bm{x})$ is also primitive. Indeed, for
every prime $p$ there is some coordinate $x_i$ not divisible by $p$, and
the coordinate $x_i^4$ of $\nu(\bm{x})$ is therefore not divisible by $p$.
Hence $\gcd(\nu(\bm{x}))=1$. It follows that $\Lambda_{\bm{x}}$ is a
primitive lattice and $\det(\Lambda_{\bm{x}})
    =
    \|\nu(\bm{x})\|$, since the determinant of a primitive lattice is equal to the determinant
of its orthogonal lattice; see, for instance,
\cite[Corollary of Lemma~1]{Schmidt1968}. The rank of $\Lambda_{\bm{x}}$ is $34$ as it is the orthogonal lattice to one fixed direction inside $\Z^{35}$.

Note that the following holds.

\begin{lemma}\label{lemma:veronese-norm}
For every $\bm{x}\in\R^4$, we have
\[
    \frac1{\sqrt{24}}\|\bm{x}\|^4
    \leq
    \|\nu(\bm{x})\|
    \leq
    \|\bm{x}\|^4.
\]
\end{lemma}

\begin{proof}
By definition,
\[
\|\nu(\bm{x})\|^2
=
\sum_{\substack{
\alpha_0+\alpha_1+\alpha_2+\alpha_3=4\\
\alpha_i\geq0
}}
x_0^{2\alpha_0}x_1^{2\alpha_1}
x_2^{2\alpha_2}x_3^{2\alpha_3}.
\]
On the other hand, by the multinomial theorem,
\begin{align*}
\|\bm{x}\|^8
=
(x_0^2+x_1^2+x_2^2+x_3^2)^4=
\sum_{|\bm{\alpha}|=4}
\binom{4}{\alpha_0,\alpha_1,\alpha_2,\alpha_3}
x_0^{2\alpha_0}x_1^{2\alpha_1}
x_2^{2\alpha_2}x_3^{2\alpha_3}.
\end{align*}
For every multi-index $\bm{\alpha}$ with $|\bm{\alpha}|=4$, we have
\[
    1
    \leq
    \binom{4}{\alpha_0,\alpha_1,\alpha_2,\alpha_3}
    \leq
    4!=24.
\]
It follows that
\[
    \|\nu(\bm{x})\|^2
    \leq
    \|\bm{x}\|^8
    \leq
    24\|\nu(\bm{x})\|^2.
\]
Taking square roots proves the claim.
\end{proof}

Let $\lambda_1(\Lambda)
    \leq
    \cdots
    \leq
    \lambda_R(\Lambda)$
denote the successive minima of a rank $R$ lattice $\Lambda$. Minkowski's second theorem (see for instance \cite[Chapter VIII, Theorem V]{Cas97}) states that
\begin{equation}
    \label{eq:Minkowskisecondtheorem}
    \det(\Lambda)\leq \lambda_1(\Lambda)
    \cdots
    \lambda_R(\Lambda)\ll \det(\Lambda).
\end{equation}

For a non-zero vector $\bm{x}\in\mathbb{Z}^4$, let
$\mathfrak{d}_2(\bm{x})$ denote the minimum determinant of a
rank-two sublattice of $\mathbb{Z}^4$ containing $\bm{x}$:
\[
\mathfrak{d}_2(\bm{x})
=
\min\left\{
\det(L):
L\subset\mathbb{Z}^4,\
\operatorname{rank}(L)=2,\
\bm{x}\in L
\right\}.
\]
We record that, for $\bm{x}\neq\bm{0}$,
\begin{equation}\label{eq:d2-trivial}
    1
    \leq
    \mathfrak{d}_2(\bm{x})
    \leq
    \|\bm{x}\|.
\end{equation}
Indeed, the lower bound follows from integrality. For the upper bound,
choose a standard basis vector $\bm{e}_i$ linearly independent from
$\bm{x}$; then the rank-two lattice generated by $\bm{x}$ and
$\bm{e}_i$ has determinant at most $\|\bm{x}\|$.

\begin{lemma}[{\cite[Lemma~3.15]{brbosaw}}]
\label{lemma:largest-minimum}
For $\bm{x}\in\Z^4_{\mathrm{prim}}$, we have
\[
    \lambda_{34}(\Lambda_{\bm{x}})
    \leq
    3
    \frac{\|\bm{x}\|}
    {\mathfrak{d}_2(\bm{x})}.
\]
\end{lemma}

We now want to revisit the proof of \cite[Lemma~3.15]{brbosaw} to get more information on some of
the smaller successive minima.
Define \[
\mathfrak{d}_3(\bm{x})
:=
\min\left\{
\det(L):
L\subset\mathbb{Z}^4,\
\operatorname{rank}(L)=3,\
\bm{x}\in L
\right\}.
\]

\begin{lemma}
    \label{lem:bounds_succ_min}
    For $\bm{x}\in\Z^4_{\mathrm{prim}}$, we have
    \[
    \lambda_{20}(\Lambda_{\bm{x}})
    \ll
    \mathfrak d_3(\bm{x}),
    \qquad
    \lambda_{30}(\Lambda_{\bm{x}})
    \ll
    \frac{\mathfrak d_2(\bm{x})}{\mathfrak d_3(\bm{x})},
    \qquad
    \lambda_{34}(\Lambda_{\bm{x}})
    \ll
    \frac{\|\bm{x}\|}{\mathfrak d_2(\bm{x})}.
    \]
\end{lemma}

\begin{proof}

Let $L_{\bm{x}}
    :=
    \bm{x}^{\perp}\cap\Z^4$, and denote its successive minima by $\eta_1(\bm{x})
    \leq
    \eta_2(\bm{x})
    \leq
    \eta_3(\bm{x}).$ We have that $L_{\bm{x}}$ is a primitive rank-three
lattice with $\det(L_{\bm{x}})=\|\bm{x}\|.$  Minkowski's second theorem \eqref{eq:Minkowskisecondtheorem} therefore gives $
    \eta_1(\bm{x})
    \eta_2(\bm{x})
    \eta_3(\bm{x})
    \asymp
    \|\bm{x}\|.$
Choose linearly independent primitive vectors $\bm{a}_1,\bm{a}_2,\bm{a}_3\in L_{\bm{x}}$ such that $\|\bm{a}_i\|
    \ll
    \eta_i(\bm{x})$, for all $i=1,2,3$. Let now
\[
    \ell_i(\bm{X})
    :=
    \bm{a}_i\cdot\bm{X}
\]
be the corresponding integral linear forms. Since
$\bm{a}_i\cdot\bm{x}=0$, we have $\ell_i(\bm{x})=0.$ Let $S=\R[X_0,X_1,X_2,X_3]$, and, for $j=1,2,3$, define
\[
    U_j
    :=
    (\ell_1,\ldots,\ell_j)_4
    \subset S_4,
\]
where the subscript denotes the degree-$4$ part of the ideal. Since the
$\ell_i$ are linearly independent, after a linear change of coordinates
we have
\[
    S/(\ell_1,\ldots,\ell_j)
    \cong
    \R[Y_1,\ldots,Y_{4-j}].
\]
It follows that
\[
    \dim U_j
    =
    35-
    \binom{4-j+3}{4}.
\]
Thus
\[
    \dim U_1=20,
    \qquad
    \dim U_2=30,
    \qquad
    \dim U_3=34.
\]
On the other hand, $U_j$ is spanned by the polynomials
\[
    m(\bm{X})\ell_i(\bm{X}),
    \qquad
    1\leq i\leq j,
\]
where $m$ runs over the monomials of degree $3$. Each such polynomial
vanishes at $\bm{x}$, and hence its coefficient vector belongs to
$\Lambda_{\bm{x}}$. Moreover, since multiplication by a monomial only
shifts the positions of the coefficients of $\ell_i$, the coefficient
vector of $m(\bm{X})\ell_i(\bm{X})$ has norm equal to that associated to $\ell_i$, which is $\|\bm{a}_i\|
    \ll
    \eta_i(\bm{x})$.
Since the above polynomials span $U_1$, we may choose $20$ linearly
independent ones whose coefficient vectors lie in $\Lambda_{\bm{x}}$ and
have norm $O(\eta_1(\bm{x}))$. Hence
\[
    \lambda_{20}(\Lambda_{\bm{x}})
    \ll
    \eta_1(\bm{x}).
\]
Similarly, the generators of $U_2$ all have coefficient-vector norm
$O(\eta_2(\bm{x}))$, since $\eta_1(\bm{x})\leq\eta_2(\bm{x})$ and
$\dim U_2=30$. Therefore
\[
    \lambda_{30}(\Lambda_{\bm{x}})
    \ll
    \eta_2(\bm{x}).
\]
Finally, using $\dim U_3=34$, we obtain
\[
    \lambda_{34}(\Lambda_{\bm{x}})
    \ll
    \eta_3(\bm{x}).
\]
Thus the construction underlying \cite[Lemma~3.15]{brbosaw} yields the
refined bounds
\begin{equation*}
    \lambda_{20}(\Lambda_{\bm{x}})
    \ll
    \eta_1(\bm{x}),
    \qquad
    \lambda_{30}(\Lambda_{\bm{x}})
    \ll
    \eta_2(\bm{x}),
    \qquad
    \lambda_{34}(\Lambda_{\bm{x}})
    \ll
    \eta_3(\bm{x})\ll\frac{\|\bm{x}\|}{\mathfrak d_2(\bm{x})} .
\end{equation*}

To finish the proof, note that we have
\[
\mathfrak d_3(\bm{x})
=
\eta_1(\bm{x}),
\qquad
\mathfrak d_2(\bm{x})
\asymp
\eta_1(\bm{x})\eta_2(\bm{x}).
\]
To see that $\mathfrak d_3(\bm{x})=\eta_1(\bm{x})$, observe that in the definition of
$\mathfrak d_3(\bm{x})$ it is enough to consider primitive rank-$3$ lattices,
since replacing a lattice by a primitive lattice containing it can only decrease its
determinant. If $\mathcal L\subset \mathbb Z^4$ is a rank-$3$ lattice
containing $\bm{x}$, then $\mathcal L^\perp$ is a primitive rank-$1$
lattice generated by some primitive vector $\bm{v}\in L_{\bm{x}}$,
and $\mathcal L=\bm{v}^{\perp}\cap\mathbb Z^4$.
Conversely, every primitive $\bm{v}\in L_{\bm{x}}$ determines in this
way a primitive rank-$3$ lattice containing $\bm{x}$. Since
$\det(\bm{v}^{\perp}\cap\mathbb Z^4)=\|\bm{v}\|$
for primitive $\bm{v}$, we obtain
\[
\mathfrak d_3(\bm{x})
=
\min_{\substack{\bm{v}\in L_{\bm{x}}\\ \bm{v}\neq 0}}
\|\bm{v}\|
=
\eta_1(\bm{x}),
\]
where a vector of minimal norm in $L_{\bm{x}}$ is automatically
primitive.
Similarly, taking orthogonal complements identifies
$\mathfrak d_2(\bm{x})$ with the minimum determinant of a primitive
rank-$2$ sublattice of $L_{\bm{x}}$. By Minkowski's second theorem
\eqref{eq:Minkowskisecondtheorem}, this minimum is comparable to the
product of the first two successive minima, and hence
\begin{equation}\label{eq:d2}
\mathfrak d_2(\bm{x})
\asymp
\eta_1(\bm{x})\eta_2(\bm{x}).
\end{equation}
This completes the proof.
\end{proof}

We will make use of the following result.

\begin{lemma}[{\cite[Lemma 3.6]{brbosaw}}]\label{lem:latticebound}
Let $N \geq 1$ and $R \in \{1,\ldots,N\}$. Let $\Lambda \subset \mathbb{R}^N$ be a lattice of rank $R$. For $T \geq 1$, we have
\[
\#\bigl(\Lambda \cap \mathcal{B}_N(T)\bigr)
\ll
\sum_{i=0}^{R}
\frac{T^{R-i}}
{\lambda_1(\Lambda)\cdots\lambda_{R-i}(\Lambda)},
\]
where the implied constant depends at most on $R$.
\end{lemma}

The last ingredient is a specialization of
\cite[Lemma~3.20]{brbosaw}, where the logarithmic factor
can be removed by keeping track of the second successive minimum in
the lattice-point count. For $X,\Delta\geq1$, define
\[
    \ell(X;\Delta)
    :=
    \#\left\{
        \bm{x}\in\Z^4:
        0<\|\bm{x}\|\leq X,\,
        \mathfrak{d}_2(\bm{x})\leq\Delta
    \right\}.
\]

\begin{lemma}\label{lemma:sparse-vectors}
For $X,\Delta\geq1$, we have
\[
    \ell(X;\Delta)
    \ll
    X^2\Delta^3.
\]
\end{lemma}

\begin{proof}
We refine the proof of \cite[Lemma~3.20]{brbosaw} in the case
$r=2$ and $n=3$. For dyadic $s_1\leq s_2$, let
$\mathcal{S}(s_1,s_2)$ denote the set of primitive rank-two
lattices $L\subset\mathbb{Z}^4$ such that
\[
    \frac{s_i}{2}<\lambda_i(L)\leq s_i,
    \qquad i=1,2.
\]
By \cite[Lemma~3.19]{brbosaw},
\begin{equation}\label{eq:lattice-count}
    \#\mathcal{S}(s_1,s_2)
    \ll
    s_1^5s_2^3.
\end{equation}
Moreover, \Cref{lem:latticebound} gives, for any rank-two lattice $L$,
\begin{equation}\label{eq:lattice-points}
\#\bigl((L\setminus\{\bm{0}\})\cap\mathcal{B}_4(X)\bigr)
\ll
\frac{X^2}{\lambda_1(L)\lambda_2(L)}
+\frac{X}{\lambda_1(L)}.
\end{equation}
Note that we can assume $\Delta\leq X^{\frac{2}{3}}$. Indeed, if $\Delta>X^{\frac{2}{3}}$ the trivial bound gives $$\ell(X;\Delta)\ll X^4
    \ll
    X^2\Delta^3.$$

For every $\bm{x}$ counted by $\ell(X;\Delta)$, choose a rank-two
lattice $\mathscr{L}_{\bm{x}}\ni\bm{x}$ with
$\det(\mathscr{L}_{\bm{x}})=\mathfrak{d}_2(\bm{x})\leq\Delta$. In particular, since $\lambda_1(\mathscr{L}_{\bm{x}})\geq1$, Minkowski's second theorem \eqref{eq:Minkowskisecondtheorem} gives $\lambda_1(\mathscr{L}_{\bm{x}})\leq \lambda_2(\mathscr{L}_{\bm{x}})\ll \Delta\leq X^{\frac{2}{3}}$. Note that $\mathscr{L}_{\bm{x}}$ is primitive, since otherwise one could find a primitive lattice containing it, which still contains
$\bm{x}$ and with smaller determinant. By Minkowski's
second theorem \eqref{eq:Minkowskisecondtheorem}, if $\mathscr{L}_{\bm{x}}\in\mathcal{S}(s_1,s_2)$, then
$s_1s_2\asymp\det(\mathscr{L}_{\bm{x}})\ll\Delta$ and $s_2\ll X$. Since $\lambda_2(\mathscr{L}_{\bm{x}})\ll X$, the second term in \eqref{eq:lattice-points} is dominated by the first one. Hence \eqref{eq:lattice-points} implies
\[
\ell(X;\Delta)
\ll
\sum_{\substack{s_1\leq s_2\\ s_1s_2\ll\Delta}}
\sum_{\mathscr{L}\in\mathcal{S}(s_1,s_2)}
\frac{X^2}{s_1s_2},
\]
where all sums over $s_1,s_2$ are dyadic.
Thus, by \eqref{eq:lattice-count}, we get
\begin{align*}
\ell(X;\Delta)
&\ll
X^2\sum_{s_2\ll X}s_2^2
\sum_{s_1\ll\min\{s_2,\frac{\Delta}{s_2}\}}s_1^4 \ll
X^2\left(
\sum_{s_2\leq\Delta^{\frac{1}{2}}}s_2^6
+
\Delta^4\sum_{s_2>\Delta^{\frac{1}{2}}}s_2^{-2}
\right)
\ll X^2\Delta^3.
\end{align*}

\end{proof}

\section{The first moment}\label{sec:first-moment}
We have
\begin{align*}
    N_F(B)
    &=
    \frac12
    \#\left\{
        \bm{x}\in\Z^4_{\mathrm{prim}}:
        0<\|\bm{x}\|\leq B,\ F(\bm{x})=0
    \right\}, \\
    \mathcal{F}_A
    &=
    \left\{
        \bm{c}\in\Z^{35}_{\mathrm{prim}}:0<\|\bm{c}\|\leq A
    \right\}/\{\pm1\}.
\end{align*}
The first moment can thus be rewritten as
\begin{align*}
\mathbb{E}_A(N_F(B))
&=
\frac{1}{2\#\mathcal{F}_A}
\sum_{\substack{
\bm{x}\in\Z^4_{\mathrm{prim}}\\
0<\|\bm{x}\|\leq B
}}
\#\left\{
F\in\mathcal{F}_A:F(\bm{x})=0
\right\}
\\
&=
\frac{1}{4\#\mathcal{F}_A}
\sum_{\substack{
\bm{x}\in\Z^4_{\mathrm{prim}}\\
0<\|\bm{x}\|\leq B
}}
\#\left\{
\bm{c}\in\Z^{35}_{\mathrm{prim}}:
\|\bm{c}\|\leq A,\,
\bm{c}\cdot\nu(\bm{x})=0
\right\}.
\end{align*}
To obtain \eqref{eq:leboudec-intro}, we therefore interpret the condition that
a quartic $F$ vanishes at a primitive vector $\bm{x}\in\Z^4$ as a lattice-point
counting problem in the primitive rank-$34$ lattice $\Lambda_{\bm{x}}$. Let $\lambda_i(\Lambda_{\bm{x}})$ be the successive minima of
$\Lambda_{\bm{x}}$. Le Boudec's combinatorial construction
\cite[Lemma~5]{leboudec2020height} gives
\[
    \lambda_{34}(\Lambda_{\bm{x}})
    \leq
    \|\bm{x}\|
    \leq B.
\]
Hence \Cref{lemma:primitive-count} applies uniformly for $A\geq B$ and
gives, as in the proof of \cite[Theorem~3]{leboudec2020height},
\begin{align*}
\mathbb{E}_A(N_F(B))
&=
\frac12
\frac{C_{34}\zeta(35)}
{C_{35}\zeta(34)}
\frac1A
\sum_{\substack{
\bm{x}\in\Z^4_{\mathrm{prim}}\\
0<\|\bm{x}\|\leq B
}}
\frac{1}{\|\nu(\bm{x})\|}
\left(
1+
O\left(\frac BA\right)
\right)
+
O\left(\frac{B^4\log(2B)}{A^{34}}\right).
\end{align*}
The next step is to analyze the weighted sum
\[
\frac{1}{A}
\sum_{\substack{\bm{x}\in\Z^4_{\mathrm{prim}}\\
0<\|\bm{x}\|\leq B}}
\frac{1}{\|\nu(\bm{x})\|}.
\]
Möbius inversion removes the primitivity condition, while Davenport's principle \cite{davenport1951lipschitz} replaces the resulting
sum by the corresponding integral within an acceptable error. Since $\nu$ is
homogeneous of degree $4$, polar coordinates in $\R^4$ separate the radial
and angular variables, and the identity $\frac{r^3}{r^4}=\frac{1}{r}$ produces a logarithmic
factor. More precisely, the angular contribution is
\[
\kappa_4
=
\int_{S^3}
\frac{d\sigma(\bm{t})}{\|\nu(\bm{t})\|},
\]
where $S^3=\partial\mathcal{B}_4(1)$ and $d\sigma(\bm{t})$ is the standard
surface measure on $S^3$ induced by the Euclidean metric on $\R^4$. Thus,
\begin{equation}\label{eq:weighted-sum}
\sum_{\substack{\bm{x}\in\Z^4_{\mathrm{prim}}\\
0<\|\bm{x}\|\leq B}}
\frac{1}{\|\nu(\bm{x})\|}
=
\frac{\kappa_4}{\zeta(4)}\log B+O(1).
\end{equation}
Combining these estimates yields \eqref{eq:leboudec-intro} (see \cite{leboudec2020height} for details).

The first observation towards improving this estimate is that we can use a
sharper bound on the largest successive minimum of
$\Lambda_{\bm{x}}$. Following \cite[(4.15)]{brbosaw}, set
\begin{equation}\label{eq:rho}
    \rho(\bm{x})
    :=
    3
    \frac{\|\bm{x}\|}
    {\mathfrak{d}_2(\bm{x})}.
\end{equation}
By \Cref{lemma:largest-minimum},
$\lambda_{34}(\Lambda_{\bm{x}})
    \leq
    \rho(\bm{x})$.
Writing the first moment as a lattice count gives
\[
    \sum_{F\in\mathcal{F}_A}N_F(B)
    =
    \frac14
    \sum_{\substack{
        \bm{x}\in\Z^4_{\mathrm{prim}}\\
        0<\|\bm{x}\|\leq B
    }}
    \mathcal{N}^*
    (\Lambda_{\bm{x}},\mathcal{B}_{35}(A)).
\]
We split the right-hand side according to whether
$\rho(\bm{x})\leq A$ or $\rho(\bm{x})>A$, and denote the two
contributions by $M^{(1)}(A,B)$ and $M^{(2)}(A,B)$, respectively:
\begin{align*}
M^{(1)}(A,B)
&:=
\frac14
\sum_{\substack{
\bm{x}\in\Z^4_{\mathrm{prim}}\\
0<\|\bm{x}\|\leq B\\
\rho(\bm{x})\leq A
}}
\mathcal{N}^*(\Lambda_{\bm{x}},\mathcal{B}_{35}(A)),
\\
M^{(2)}(A,B)
&:=
\frac14
\sum_{\substack{
\bm{x}\in\Z^4_{\mathrm{prim}}\\
0<\|\bm{x}\|\leq B\\
\rho(\bm{x})>A
}}
\mathcal{N}^*(\Lambda_{\bm{x}},\mathcal{B}_{35}(A)).
\end{align*}
The first sum will be referred to as the \emph{good contribution},
while the second one will be referred to as the \emph{bad contribution}.

\subsection*{The good contribution}

Suppose that $\rho(\bm{x})\leq A$. We may apply
\Cref{lemma:primitive-count} with $Y=\rho(\bm{x})$ to obtain
\begin{align*}
    \mathcal{N}^*
    (\Lambda_{\bm{x}},\mathcal{B}_{35}(A))
    =
    \frac{C_{34}}{\zeta(34)}
    \frac{A^{34}}{\|\nu(\bm{x})\|}
    +
    O\left(
        A^{33}
        \frac{\|\bm{x}\|}
        {\mathfrak{d}_2(\bm{x})\|\nu(\bm{x})\|}
    \right)
    +
    O(A\log(2B)).
\end{align*}
Hence
\begin{align}\label{eq:good-contribution}
    M^{(1)}(A,B)
    =&
    \frac{C_{34}}{4\zeta(34)}
    A^{34}
    \sum_{\substack{
        \bm{x}\in\Z^4_{\mathrm{prim}}\\
        0<\|\bm{x}\|\leq B\\
        \rho(\bm{x})\leq A
    }}
    \frac1{\|\nu(\bm{x})\|}
    +
    O\left(A^{33}T(B)\right)
    +
    O\left(AB^4\log(2B)\right),
\end{align}
where
\[
    T(B)
    :=
    \sum_{\substack{
        \bm{x}\in\Z^4_{\mathrm{prim}}\\
        0<\|\bm{x}\|\leq B
    }}
    \frac{\|\bm{x}\|}
    {\mathfrak{d}_2(\bm{x})\|\nu(\bm{x})\|}.
\]
Since we are inside an error term, we have harmlessly enlarged the sum by
dropping the condition $\rho(\bm{x})\leq A$.

We estimate $T(B)$ dyadically. On a shell
\[
    X<\|\bm{x}\|\leq2X,
\]
split
\[
    D<\mathfrak{d}_2(\bm{x})\leq2D.
\]
Recall that by \eqref{eq:d2-trivial} we have
$\mathfrak{d}_2(\bm{x})\leq\|\bm{x}\|$, so that it is enough to consider
$D\ll X$. By \Cref{lemma:sparse-vectors} and the trivial estimate
$O(X^4)$, the number of vectors in such a box is
\[
    \#\left\{
        \bm{x}\in\Z^4:
        X<\|\bm{x}\|\leq2X,\,
        D<\mathfrak{d}_2(\bm{x})\leq2D
    \right\}
    \ll
    \min\{X^4,X^2D^3\}.
\]
Consequently,
\begin{align*}
\sum_{\substack{
X<\|\bm{x}\|\leq2X\\
\bm{x}\in\Z^4_{\mathrm{prim}}
}}
\frac1{\mathfrak{d}_2(\bm{x})}
&\ll
\sum_{\substack{D\ll X\\D\text{ dyadic}}}
\frac1D
\min\left\{
X^4,
X^2D^3
\right\}.
\end{align*}
We split the last sum at $D_0=X^{\frac{2}{3}}$. For $D\leq D_0$, we have
\begin{align*}
\sum_{\substack{
D\leq D_0\\
D\text{ dyadic}
}}
\frac1D
X^2D^3
&\ll
X^2
\sum_{\substack{
D\leq D_0\\
D\text{ dyadic}
}}
D^2
\ll
X^2D_0^2
\ll
X^{\frac{10}{3}}.
\end{align*}
For $D>D_0$,
\begin{align*}
\sum_{\substack{
D>D_0\\
D\ll X\\
D\text{ dyadic}
}}
\frac{X^4}{D}
&\ll
\frac{X^4}{D_0}
\ll
X^{\frac{10}{3}}.
\end{align*}
We have therefore proved that
\begin{equation}\label{eq:d2-average-shell}
\sum_{\substack{
X<\|\bm{x}\|\leq2X\\
\bm{x}\in\Z^4_{\mathrm{prim}}
}}
\frac1{\mathfrak{d}_2(\bm{x})}
\ll
X^{\frac{10}{3}}.
\end{equation}
Using \Cref{lemma:veronese-norm}, on the same dyadic shell we have
\[
    \frac{\|\bm{x}\|}{\|\nu(\bm{x})\|}
    \ll
    \frac1{X^3}.
\]
Thus \eqref{eq:d2-average-shell} yields
\begin{align*}
\sum_{\substack{
X<\|\bm{x}\|\leq2X\\
\bm{x}\in\Z^4_{\mathrm{prim}}
}}
\frac{\|\bm{x}\|}
{\mathfrak{d}_2(\bm{x})\|\nu(\bm{x})\|}
&\ll
X^{-3}
X^{\frac{10}{3}}
\ll
X^{\frac{1}{3}}.
\end{align*}
Summing over the dyadic intervals $X=2^k\ll B$, we obtain
\begin{equation}\label{eq:T-bound}
    T(B)\ll B^{\frac{1}{3}}.
\end{equation}
Hence the total error coming from the successive-minima term in the
good range is $O\left(A^{33}B^{\frac{1}{3}}\right)$. After dividing by $\#\mathcal{F}_A\asymp A^{35}$, this produces the normalized
error $O\left(
        \frac{B^{\frac{1}{3}}}{A^2}
    \right)$. The additional term $O(AB^4\log(2B))$ in
\eqref{eq:good-contribution} contributes $O\left(\frac{B^4\log(2B)}{A^{34}}\right)$, which will be absorbed by the later error terms.

\subsection*{The bad contribution}

We now assume that $\rho(\bm{x})>A$. By \eqref{eq:rho} we have
\begin{equation}\label{eq:bad-d2}
    \mathfrak{d}_2(\bm{x})
    \ll
    \frac{\|\bm{x}\|}{A}.
\end{equation}
We divide the vectors into dyadic boxes
\[
    X<\|\bm{x}\|\leq2X,
    \qquad
    D<\mathfrak{d}_2(\bm{x})\leq2D.
\]
Then $D\ll\frac{X}{A}$, and, since $\mathfrak d_2(\bm{x})\geq1$, there are no bad boxes unless
$X\gg A$. Since
\[
    \mathcal{N}^*(\Lambda_{\bm{x}},\mathcal{B}_{35}(A))
    \leq
    \#\bigl(\Lambda_{\bm{x}}\cap\mathcal{B}_{35}(A)\bigr),
\]
\Cref{lem:latticebound}, with $R=34$ and $T=A$, gives after reindexing
\begin{equation}\label{eq:lemmabound}
    \mathcal{N}^*(\Lambda_{\bm{x}},\mathcal{B}_{35}(A))
    \ll
    \sum_{i=0}^{34}
    \frac{A^i}
    {\lambda_1(\Lambda_{\bm{x}})\cdots\lambda_i(\Lambda_{\bm{x}})},
\end{equation}
where the denominator is interpreted as $1$ for $i=0$. On the $(X,D)$-box, \Cref{lem:bounds_succ_min} and
$\mathfrak d_3(\bm{x})\geq1$ give
\[
    \lambda_{30}(\Lambda_{\bm{x}})\ll D,
    \qquad
    \lambda_{34}(\Lambda_{\bm{x}})\ll\frac{X}{D}.
\]
Moreover, by Minkowski's second theorem \eqref{eq:Minkowskisecondtheorem} and
\Cref{lemma:veronese-norm}, $\lambda_1(\Lambda_{\bm{x}})
    \cdots
    \lambda_{34}(\Lambda_{\bm{x}})
    \asymp X^4.$ Hence, for every $i\in\{0,\ldots,34\}$,
\begin{equation}\label{eq:bad-complementary-minima}
    \frac{A^i}
    {\lambda_1(\Lambda_{\bm{x}})\cdots\lambda_i(\Lambda_{\bm{x}})}
    \ll
    \frac{A^i}{X^4}
    \lambda_{i+1}(\Lambda_{\bm{x}})
    \cdots
    \lambda_{34}(\Lambda_{\bm{x}}).
\end{equation}

We estimate the terms in \eqref{eq:lemmabound} in four ranges. If
$30\leq i\leq34$, then
\[
    \lambda_{i+1}(\Lambda_{\bm{x}})
    \cdots
    \lambda_{34}(\Lambda_{\bm{x}})
    \ll
    \left(\frac{X}{D}\right)^{34-i},
\]
and therefore
\[
    \frac{A^i}
    {\lambda_1(\Lambda_{\bm{x}})\cdots\lambda_i(\Lambda_{\bm{x}})}
    \ll
    \frac{A^{30}}{D^4}
    \left(\frac{AD}{X}\right)^{i-30}
    \ll
    \frac{A^{30}}{D^4},
\]
since $AD\ll X$.
If $28\leq i\leq29$, monotonicity of the successive minima gives
\[
    \lambda_{i+1}(\Lambda_{\bm{x}})
    \cdots
    \lambda_{30}(\Lambda_{\bm{x}})
    \ll
    D^{30-i},
\]
while
\[
    \lambda_{31}(\Lambda_{\bm{x}})
    \cdots
    \lambda_{34}(\Lambda_{\bm{x}})
    \ll
    \left(\frac{X}{D}\right)^4.
\]
Thus
\[
    \frac{A^i}
    {\lambda_1(\Lambda_{\bm{x}})\cdots\lambda_i(\Lambda_{\bm{x}})}
    \ll
    A^iD^{26-i}.
\]
If $26\leq i\leq27$, we keep track of $\mathfrak d_3(\bm{x})$ instead of using
$\mathfrak d_3(\bm{x})\geq1$. By \Cref{lem:bounds_succ_min} and monotonicity of the
successive minima, we have
\[
    \lambda_{i+1}(\Lambda_{\bm{x}})
    \cdots
    \lambda_{30}(\Lambda_{\bm{x}})
    \ll
    \left(\frac{\mathfrak d_2(\bm{x})}{\mathfrak d_3(\bm{x})}\right)^{30-i},
    \qquad
    \lambda_{31}(\Lambda_{\bm{x}})
    \cdots
    \lambda_{34}(\Lambda_{\bm{x}})
    \ll
    \left(\frac{X}{\mathfrak d_2(\bm{x})}\right)^4,
\]
and therefore \eqref{eq:bad-complementary-minima} gives
\[
    \frac{A^{i}}
    {\lambda_1(\Lambda_{\bm{x}})\cdots\lambda_{i}(\Lambda_{\bm{x}})}
    \ll
    \frac{A^{i}}{\mathfrak d_3(\bm{x})^{30-i}\mathfrak d_2(\bm{x})^{i-26}}
    \ll
    \frac{A^{27}}{\mathfrak d_3(\bm{x})^4},
\]
where in the last step we have used the fact that
$\mathfrak d_2(\bm{x})\gg\mathfrak d_3(\bm{x})$, which follows from
\eqref{eq:d2} since $\mathfrak d_3(\bm{x})=\eta_1(\bm{x})$ and $\eta_2(\bm{x})\geq1$.
Finally, for $i\leq25$, integrality gives
$\lambda_j(\Lambda_{\bm{x}})\geq1$, and hence
\[
    \frac{A^i}
    {\lambda_1(\Lambda_{\bm{x}})\cdots\lambda_i(\Lambda_{\bm{x}})}
    \ll
    A^{25}.
\]
Combining these estimates with \eqref{eq:lemmabound}, we obtain
\[
    \mathcal{N}^*(\Lambda_{\bm{x}},\mathcal{B}_{35}(A))
    \ll
    A^{25}
    +
    \frac{A^{27}}{\mathfrak d_3(\bm{x})^4}
    +
    \frac{A^{28}}{D^2}
    +
    \frac{A^{29}}{D^3}
    +
    \frac{A^{30}}{D^4}.
\]
Using \Cref{lemma:sparse-vectors} for all the terms except the second one, it follows that
\begin{align*}
M^{(2)}(A,B)
\ll{}&
A^{27}
\sum_{X\ll B}
\sum_{\substack{
X<\|\bm{x}\|\leq2X\\
\bm{x}\in\Z^4_{\mathrm{prim}}
}}
\frac{1}{\mathfrak d_3(\bm{x})^4}
\\
&+
\sum_{\substack{X\ll B\\ D\ll \frac{X}{A}}}
\biggl(
    X^2D^3A^{25}
    +X^2DA^{28}
    +X^2A^{29}
    +\frac{X^2A^{30}}{D}
\biggr),
\end{align*}
where all sums over $X$ and $D$ are dyadic, and where in the first sum we have
harmlessly dropped the condition $\rho(\bm{x})>A$. To estimate the first
sum, we split $E\leq\mathfrak d_3(\bm{x})<2E$ for dyadic $E\geq1$ and apply
\cite[Lemma~3.20]{brbosaw} in the case $r=3$ and $n=3$, which shows that the
number of non-zero $\bm{x}\in\Z^4$ with $\|\bm{x}\|\leq2X$ and
$\mathfrak d_3(\bm{x})\leq2E$ is $O(X^3E^3\log(2E))$. Here the logarithmic
factor is harmless, and we obtain
\[
\sum_{\substack{
X<\|\bm{x}\|\leq2X\\
\bm{x}\in\Z^4_{\mathrm{prim}}
}}
\frac{1}{\mathfrak d_3(\bm{x})^4}
\ll
\sum_{E\text{ dyadic}}
\frac{1}{E^4}X^3E^3\log(2E)
\ll
X^3.
\]
Summing over $D$ in the second sum, we therefore obtain
\[
    M^{(2)}(A,B)
    \ll
    \sum_{X\ll B}
    \left(
        X^5A^{22}
        +X^3A^{27}
        +X^2A^{29}\log(X)
        +X^2A^{30}
    \right).
\]
Summing over the dyadic values of $X$ and using
$\#\mathcal F_A\asymp A^{35}$, we conclude that
\begin{equation}\label{eq:M2-normalized}
    \frac{M^{(2)}(A,B)}{\#\mathcal F_A}
    \ll
    \frac{B^5}{A^{13}}
    +
    \frac{B^3}{A^8}
    +
    \frac{B^2\log(B)}{A^6}
    +
    \frac{B^2}{A^5}.
\end{equation}

\begin{remark}
The estimate \eqref{eq:M2-normalized} is not optimal: in the range
$i\leq25$ we have only used integrality, and the bounds of
\Cref{lem:bounds_succ_min} could be exploited further to improve the term
$\frac{B^5}{A^{13}}$. However, this term only exceeds $\frac{B^2}{A^5}$ when $B>A^{\frac{8}{3}}$,
in which region the latter term is already larger than one.
Such a refinement would therefore benefit neither \Cref{cor:first-moment}
nor \Cref{thm:smallest-point}.
\end{remark}

\subsection*{Completion of the proof}

The main term in \eqref{eq:good-contribution} only contains vectors with
$\rho(\bm{x})\leq A$. To restore the omitted vectors, note from
\eqref{eq:bad-d2}, \Cref{lemma:sparse-vectors} and
\Cref{lemma:veronese-norm} that, on a dyadic shell
$X<\|\bm{x}\|\leq2X,$
we have
\[
    \sum_{\substack{
        X<\|\bm{x}\|\leq2X\\
        \bm{x}\in\Z^4_{\mathrm{prim}}\\
        \rho(\bm{x})>A
    }}
    \frac1{\|\nu(\bm{x})\|}
    \ll
    \frac1{X^4}
    \ell\left(2X;\frac{6X}{A}\right)
    \ll
    \frac{X}{A^3}.
\]
Note that, by \eqref{eq:d2-trivial}, the sum on the left-hand side is empty unless $\frac{6X}{A}\geq1$, so that \Cref{lemma:sparse-vectors} is indeed applicable.
Therefore
\begin{equation}\label{eq:missing-main}
    \sum_{\substack{
        \bm{x}\in\Z^4_{\mathrm{prim}}\\
        0<\|\bm{x}\|\leq B\\
        \rho(\bm{x})>A
    }}
    \frac1{\|\nu(\bm{x})\|}
    \ll
    \frac{B}{A^3}.
\end{equation}
After normalization by $\frac{1}{A}$ this contributes $O\left(\frac{B}{A^4}\right)$.
Combining
\eqref{eq:good-contribution},
\eqref{eq:T-bound},
\eqref{eq:M2-normalized},
\eqref{eq:missing-main},
\eqref{eq:weighted-sum}
and
\eqref{eq:family-size},
we obtain
\[ \mathbb{E}_A(N_F(B)) = c_4\frac{\log B}{A} + O\left( \frac1A + \frac{B^{\frac{1}{3}}}{A^2} + \frac{B}{A^4} + \frac{B^2}{A^5} + \frac{B^2\log B}{A^6} + \frac{B^3}{A^8} + \frac{B^5}{A^{13}} \right). \] The additional error $O\left(\frac{B^4\log(2B)}{A^{34}}\right) $ coming from the $O(A\log(2B))$ term in the good contribution is absorbed by $O(\frac{B^5}{A^{13}})$. Similarly, replacing $\frac{1}{\#\mathcal{F}_A}$ in the main term by its main term from \eqref{eq:family-size} introduces an error of $O\left( \frac{\log B}{A^2}\right)$ which is absorbed by $O(\frac{B^{\frac{1}{3}}}{A^2})$. Moreover, observe that for $B\leq A^3$ we have $\frac{B^{1/3}}{A^2}\leq \frac{(A^3)^{1/3}}{A^2}=\frac{1}{A}$, while for $B>A^3$ we have $\frac{B^{1/3}}{A^2}=\frac{B^{1/3}}{A^2}\cdot \frac{A^3}{A^3}\leq \frac{B^{4/3}}{A^5}\leq \frac{B^2}{A^5}$. All in all, we get \[
\frac{B^{1/3}}{A^2}
\ll \frac1A+\frac{B^2}{A^5}.
\]
Similarly, splitting at $B=A$, we have \[
\frac{B}{A^4}
\ll \frac1A+\frac{B^2}{A^5}.
\]
On the other hand, splitting at $\log B=A$, we have
\[
\frac{B^2\log B}{A^6}
\ll
\frac{B^2}{A^5}+\frac{B^3}{A^8},
\]
and also \[
\frac{B^3}{A^8}
\ll
\frac{B^2}{A^5}+\frac{B^5}{A^{13}},
\]
splitting at $B=A^{\frac{5}{2}}$. This proves \Cref{thm:first-moment}. To deduce \Cref{cor:first-moment}, we compare the error terms with the main term $\frac{\log B}{A}$. The condition $ B^2=o(A^4\log B) $ is precisely what is required for the term $\frac{B^2}{A^5}$, since \[ \frac{\frac{B^2}{A^5}}{\frac{\log B}{A}} = \frac{B^2}{A^4\log B}. \] Moreover, this condition implies $ B\ll A^2(\log B)^{\frac{1}{2}}$, so that $\log B\ll\log A$, and hence  \[ \frac{\frac{B^5}{A^{13}}}{\frac{\log B}{A}} = \frac{B^2}{A^4\log B} \frac{B^3}{A^8} = o(1). \] Also, \[ \frac{\frac1A}{\frac{\log B}{A}} = \frac1{\log B} = o(1). \] Therefore, \Cref{cor:first-moment} follows.

\section{Minimal height of a rational point}
\label{sec:minimalheight}
In this section, we follow the strategy of \cite[Theorem~1]{leboudec2020height}
to obtain information on the smallest height of a rational point on surfaces
in $\mathbb{V}$. Our first-moment estimate can be used to show that quartic surfaces admitting a rational point of height
$o(H(V)^{\frac{5}{2}})$ form a density-zero subset of the family.

\begin{proof}[Proof of \Cref{thm:smallest-point}]
It is enough to prove that
\[
\frac{
\#\{V\in\mathcal F_A:M(V)\leq\psi(H(V))\}
}{
\#\mathcal F_A
}
\longrightarrow 0.
\]
Fix $\varepsilon\in(0,1)$. Since
$\psi(u)=o(u^{\frac{5}{2}})$, for $A$ sufficiently large we have $\psi(u)\leq \varepsilon u^{\frac{5}{2}}$
for every $u\geq A^{\frac{1}{2}}$. Hence, for
$A^{\frac{1}{2}}<H(V)\leq A$,
\[
    \psi(H(V))
    \leq
    \varepsilon H(V)^{\frac{5}{2}}
    \leq
    \varepsilon A^{\frac{5}{2}},
\]
and therefore
\begin{align*}
\#\{V\in\mathcal F_A:M(V)\leq\psi(H(V))\}
\leq
\#\mathcal F_{A^{\frac{1}{2}}}
+
\#\left\{
V\in\mathcal F_A:
M(V)\leq \varepsilon A^{\frac{5}{2}}
\right\}.
\end{align*}

If $M(V_F)\leq B$, then $N_F(B)\geq1$. Thus
\[
\#\left\{
V_F\in\mathcal F_A:
M(V_F)\leq \varepsilon A^{\frac{5}{2}}
\right\}
\leq
\sum_{V_F\in\mathcal F_A}
N_F\left(\varepsilon A^{\frac{5}{2}}\right).
\]
Dividing by $\#\mathcal F_A$ and using
\eqref{eq:family-size}, we obtain
\[
\frac{
\#\{V\in\mathcal F_A:M(V)\leq\psi(H(V))\}
}{
\#\mathcal F_A
}
\ll
A^{-\frac{35}{2}}
+
\mathbb E_A\left(
N_F\left(\varepsilon A^{\frac{5}{2}}\right)
\right).
\]
We now apply \Cref{thm:first-moment} with $B=\varepsilon A^{\frac{5}{2}}$.
This gives
\begin{align*}
\mathbb E_A(N_F(B))
\ll{}&
\frac{\log A}{A}
+\varepsilon^2
+\frac{\varepsilon^5}{A^{\frac{1}{2}}}.
\end{align*}
For fixed $\varepsilon$, the first and the last term on the right-hand side
tend to zero as $A\to\infty$. We deduce that
\[
\limsup_{A\to\infty}
\mathbb E_A\left(
N_F\left(\varepsilon A^{\frac{5}{2}}\right)
\right)
\ll
\varepsilon^2;\ \text{ and } \
\limsup_{A\to\infty}
\frac{
\#\{V\in\mathcal F_A:M(V)\leq\psi(H(V))\}
}{
\#\mathcal F_A
}
\ll
\varepsilon^2.
\]
Since $\varepsilon>0$ is arbitrary, the desired limit follows.
\end{proof}

\printbibliography

\end{document}